\documentclass[11pt]{article}
\usepackage{amsmath, amsfonts, amssymb, amsthm}
\usepackage{natbib}
\usepackage{hyperref}

\DeclareMathSymbol{\leqslant}{\mathalpha}{AMSa}{"36} 

\DeclareMathSymbol{\geqslant}{\mathalpha}{AMSa}{"3E} 

\DeclareMathSymbol{\eset}{\mathalpha}{AMSb}{"3F}     

\newcommand{\suptwo}[2]{\sup_{\substack{#1 \\ #2}}} 

\newcommand{\range}{\mathrm{\,range\,}}

\newcommand{\cA}{\mathcal{A}}

\newcommand{\cC}{\mathcal{C}}

\newcommand{\cF}{\mathcal{F}}

\newcommand{\cL}{\mathcal{L}}

\newcommand{\cO}{\mathcal{O}}

\newcommand{\cT}{\mathcal{T}}

\newcommand{\Om}{\Omega}

\newcommand{\eps}{\varepsilon}
\newcommand{\lam}{\lambda}

\newcommand{\vphi}{\varphi}

\newcommand{\tk}{\tilde\kappa}

\newcommand{\lfl}{\lfloor}
\newcommand{\rfl}{\rfloor}

\newcommand{\dto}{\downarrow}

\newcommand{\im}{\mathrm{im}}

\newcommand{\fB}{X}
\newcommand{\hB}{\hat X}

\newcommand{\IP}{\mathbb{P}}

\newcommand{\IL}{\mathbb{D}}
\newcommand{\ID}{\mathbb{D}}
\newcommand{\IN}{\mathbb{N}}
\newcommand{\IH}{\mathbb{H}}
\newcommand{\IZ}{\mathbb{Z}}

\newcommand{\IC}{\mathbb{C}}
\newcommand{\IR}{\mathbb{R}}
\newcommand{\IE}{\mathbb{E}}
\newcommand{\iN}{\in\IN}
\newcommand{\iZ}{\in\IZ}
\newcommand{\iR}{\in\IR}

\newcommand{\be}{\begin{eqnarray*}}
\newcommand{\ee}{\end{eqnarray*}}
\newcommand{\ben}{\begin{eqnarray}}
\newcommand{\een}{\end{eqnarray}}

\theoremstyle{plain}
\newtheorem{theo}{Theorem}[section]
\newtheorem{lemma}[theo]{Lemma}
\newtheorem{propo}[theo]{Proposition}

\theoremstyle{definition}

\newtheorem{remark}[theo]{Remark}

\renewenvironment{proof}[1][] {\noindent{\bf Proof#1. }}{\hspace*{\fill}$\square$\medskip\par}

\begin{document}

\vglue20pt
\centerline{\huge\bf High-resolution }
\medskip

\centerline{\huge\bf quantization and entropy coding} 
\medskip

\centerline{\huge\bf for fractional Brownian motion}

\bigskip
\bigskip

\centerline{by}
\bigskip
\medskip

\centerline{{\Large S. Dereich and  M. Scheutzow}}

\bigskip

\centerline{\it  Technische Universit\"at Berlin}

\bigskip
\bigskip
\bigskip

{\leftskip=1truecm
\rightskip=1truecm
\baselineskip=15pt
\small

\noindent{\slshape\bfseries Summary.} 
We derive a high-resolution formula for the quantization and entropy coding approximation quantities for fractional Brownian motion, respective to the supremum norm and $L^p[0,1]$-norm distortions. We show that all moments in the quantization problem lead to the same asymptotics. Using a general principle, we conclude that entropy coding and quantization coincide asymptotically.
Under supremum-norm distortion, our proof uses an explicit construction of efficient codebooks based on a particular entropy constrained coding scheme.
This procedure can be used to construct close to optimal high resolution quantizers.

\bigskip


\bigskip


}

\section{Introduction}

Functional quantization and entropy coding  concern the finding of ``good" discrete approximations to a non-discrete random signal in a Banach space of functions. 
Such discrete approximations may serve as evaluation points for quasi Monte Carlo methods or as an information reduction of the original to allow storage on a computer or transmission over some channel with finite capacity. In the past years, research in this field has been very active, which resulted in  numerous new results. 
Previous research addressed, for instance, the problem of constructing good approximation schemes, the evaluation of the theoretically best approximation under an information constraint, existence of optimal approximation schemes and regularity properties of the paths of optimal approximations. 
The above questions are treated for Gaussian measures in Hilbert spaces  by Luschgy and Pag\`es (\cite{LuPa02}, \cite{LuPa04}) and by the first-named author in \cite{Der03}. For Gaussian originals in Banach spaces, these problems have been addressed by the authors and collaborators in \cite{DFMS03}, \cite{DerLif04}, \cite{Der03}, \cite{Der04} and by Graf, Luschgy and Pag\`es in \cite{GLP05}.
For  general accounts of quantization and coding theory in finite dimensional spaces, see \cite{GraLu00} and \cite{CoTho91} (see also \cite{GrNeu98}).

In this article, we consider the asymptotic coding problem of fractional Brownian motion for the supremum and $L^p[0,1]$-norm distortions. We derive the asymptotic quality of  optimal approximations. In particular, it is shown that efficient entropy constrained quantizers can be used to construct close to optimal quantizers when considering the supremum norm. Moreover, for one of the above norm-based distortions, all moments and both information constraints lead to the same asymptotic approximation quality. In particular, quantization is asymptotically just as efficient as entropy coding. 
The main impetus to the present work was provided by the necessity to understand the coding complexity of Brownian motion in order to solve the quantization (resp. entropy constrained coding) problem for diffusions (see \cite{Der04b}).

%
Let $(\Om,\cA,\IP)$ be a probability space, let $H\in(0,1)$ and let $\fB=(\fB_t)_{t\ge0}$ denote fractional Brownian motion with Hurst index $H$ on $(\Om,\cA,\IP)$, i.e.\ $(\fB_t)_{t\ge0}$ is a centered continuous Gaussian process with covariance kernel
$$
K(t,s)=\frac 12 [t^{2H}+s^{2H}-|t-s|^{2H}],\qquad t,s\ge0.
$$
We need some more notation. 
In the sequel, $\IC[0,a]$, $a>0$, and $\IL[0,a]$ denote the space of continuous real-valued functions on the interval $[0,a]$ and the space of c\`adl\`ag functions on $[0,a]$, respectively. Both spaces are endowed with the supremum norm $\|\cdot\|_{[0,a]}$.
Moreover, we let $(L^p[0,a],\|\cdot\|_{L^p[0,a]})$ denote the standard $L^p$-space of real-valued functions defined on $[0,a]$. 
Finally, $\|\cdot\|_q$, $q\in(0,\infty]$, denotes the  $L^q$-norm induced by the  probability measure $\IP$ on the set of real-valued random variables.

Let us briefly introduce the main objectives of quantization and entropy coding.
Let $E$ and $\hat E$ denote measurable spaces, and let $d:E\times \hat E\to [0,\infty)$ be a product measurable function. 
For a given $E$-valued r.v.\ $Y$ (\emph{original}) and \emph{moment} $q>0$, the aim is to minimize
\begin{align}\label{eq0214-1}
\bigl\|d(Y,\pi(Y))\bigr\|_q 
\end{align}
over all measurable functions $\pi:E\to \hat E$ with discrete image (\emph{strategy}) that satisfy a particular information constraint parameterized by the \emph{rate} $r\ge0$.

\emph{Entropy coding} (also known as \emph{entropy constrained quantization} in the literature) concerns the minimization of (\ref{eq0214-1}) over all strategies $\pi$ having  entropy $\IH(\pi(Y))$ at most $r$. Recall that the entropy of a discrete r.v.\ $Z$ with probability weights $(p_w)$ is defined as
$$
\IH(Z)=-\sum_w p_w \log p_w=\IE[-\log p_Z].
$$
In the \emph{quantization problem}, one is considering strategies $\pi$ satisfying the range constraint: $|\range(\pi(Y))|\le e^r$.
The corresponding approximation quantities are the \emph{entropy-constrained quantization error}
\begin{align}\label{eq0217-1}
D^{(e)}(r|Y,E,\hat E,d,q):=\inf_\pi \ \bigl\| d(Y,\pi(Y)) \bigr\|_q, 
\end{align}
where the infimum is taken over all strategies $\pi$ with entropy rate $r\ge0$, and
the \emph{quantization error}
\begin{align}\label{eq0217-2}
D^{(q)}(r|Y,E, \hat E,d,q):=\inf_\pi \ \bigl\| d(Y,\pi(Y)) \bigr\|_q, 
\end{align}
the infimum being taken over all strategies $\pi$ having quantization rate $r\ge0$. Often, all or some of the parameters $Y$, $E$, $\hat E$, $d$, $q$ are clear from the context. 
Then we omit these parameters in the quantities $D^{(e)}$ and $D^{(q)}$.
The \emph{quantization information constraint} is more restrictive, so that the quantization error always dominates the entropy coding error. Moreover, the coding error increases with the moment under consideration. 

Unless otherwise stated, we choose as original $Y=X$ and as original space $E=\IC[0,\infty)$.
We are mainly concerned with two particular choices for $\hat E$ and $d$. In the first sections, we treat the case where $\hat E=\ID[0,1]$ and $d(f,g)=\|f-g\|_{[0,1]}$. In this setting we find:

\begin{theo} \label{theo0303-1}
There exists a constant $\kappa=\kappa(H)\in(0,\infty)$ such that for all $q_1\in(0,\infty]$ and $q_2\in(0,\infty)$,
$$
\lim_{r\to\infty} r^H \, D^{(e)}(r|q_1)= \lim_{r\to\infty} r^H \, D^{(q)}(r|q_2) = \kappa.
$$
\end{theo}

\begin{remark} In the above theorem, general c\`adl\`ag functions are allowed as reconstructions. Since the original process is continuous, it might seem more natural to use continuous functions as approximations. The following argument shows that, for a finite moment $q>0$, the space $\hat E=\ID[0,1]$ can be replaced by $\hat E=\IC[0,1]$ without changing $D^{(q)}$ and $D^{(e)}$. Let $\pi:\IC[0,1]\to\ID[0,1]$ be an arbitrary strategy and let $\tau_n:\ID[0,1]\to\IC[0,1]$ denote the linear operator mapping $f$ to its piecewise linear interpolation with supporting points $0,\frac1n,\frac2n\dots,1$. Then
\begin{align*}
\bigl\|\|X-\tau_n\circ\pi(X)\|_{[0,1]} \bigr\|_q &\le \bigl\|\|\tau_n (X)-\tau_n\circ\pi(X)\|_{[0,1]} \bigr\|_q+\bigl\|\|X- \tau_n (X)\|_{[0,1]} \bigr\|_q\\
&\le \bigl\|\|X-\pi(X)\|_{[0,1]} \bigr\|_q+\bigl\|\|X- \tau_n (X)\|_{[0,1]} \bigr\|_q.
\end{align*}
Note that the second term vanishes when $n$ tends to infinity and that $\tau_n\circ \pi$ satisfies the same information constraint as $\pi$.  
\end{remark}

In the last section we conclude the article with a discussion of the case where $\hat E=L^p[0,1]$ and $d(f,g)=\|f-g\|_{L^p[0,1]}$ for some $p\ge1$. In this case, one has the following analog to Theorem \ref{theo0303-1}:

\begin{theo} \label{theo0303-2}
For every $p\ge1$ there exists a constant $\kappa=\kappa(H,p)\in(0,\infty)$ such that for all $q\in(0,\infty)$,
$$
\lim_{r\to\infty} r^H \, D^{(e)}(r|q)= \lim_{r\to\infty} r^H \, D^{(q)}(r|q) = \kappa.
$$
\end{theo}

\begin{remark}
It is again possible to replace the space $\hat E=L^p[0,1]$ by $\hat E=\IC[0,1]$ without changing $D^{(q)}$ and $D^{(e)}$. Indeed, for $\eps>0$, let  $h_\eps:\IR\to[0,\infty)$ denote a smooth function supported on $[-\eps,\eps]$ with $\int f_\eps=1$, and define $\tau_\eps:L^p[0,1]\to\IC[0,1]$ through $\tau_\eps(f)(t)=\int_0^1 f(s) \, h(t-s) \,ds$. Then for a given strategy $\pi:\IC[0,1]\to L^p[0,1]$ one obtains 
\begin{align*}
\bigl\|\|X-\tau_\eps\circ\pi(X)\|_{L^p[0,1]} \bigr\|_q &\le \bigl\|\|\tau_\eps (X)-\tau_\eps\circ\pi(X)\|_{L^p[0,1]} \bigr\|_q+\bigl\|\|X- \tau_\eps (X)\|_{L^p[0,1]} \bigr\|_q\\
&\le \bigl\|\|X-\pi(X)\|_{L^p[0,1]} \bigr\|_q+\bigl\|\|X- \tau_\eps (X)\|_{L^p[0,1]} \bigr\|_q,
\end{align*}
where the last inequality is a consequence of Young's inequality. Now for $\eps\dto0$ the second term converges to $0$.
\end{remark}

For ease of notation, the article is restricted to the analysis of $1$-dimensional processes. 
However, when replacing $(X_t)$ by a process $(X^{(1)}_t,\dots,X^{(d)}_t)$ consisting of $d$ independent fractional Brownian motions, the proofs can be easily adapted, and one obtains analogous results. In particular, it is possible to prove analogs of the above theorems for a multi dimensional Brownian motion.

Let us summarize some of the known estimates for the constant $\kappa$ in the case where $X$ is standard Brownian motion, i.e. $H=1/2$.
\begin{itemize}
\item When $\hat E=\ID[0,1]$ and $d(f,g)=\|f-g\|_{[0,1]}$, the relationship between the small ball function and the quantization problem (see \cite{DFMS03}) leads to
$$
\kappa\in \bigl[\frac\pi{\sqrt 8}, \pi\bigr].
$$ 
\item For $\hat E=L^p[0,1]$, $p\ge1$, and $d(f,g)=\|f-g\|_{[L^p[0,1]}$, $\kappa$ may again be estimated via a connection to the small ball function. Indeed, letting
$$
\lam_1=\inf\Bigl\{\int_{-\infty}^\infty |x|^p \vphi^2(x)\,dx+{\textstyle\frac12} \int_{-\infty}^\infty (\vphi'(x))^2\,dx\Bigr\},
$$
where the infimum is taken over all weakly differentiable $\vphi\in L^2(\IR)$ with unit norm, one has
$$
\kappa\in[c,\sqrt 8 \,c]
$$
for $c=2^{1/p} \sqrt p\, \bigl(\frac{\lam_1}{2+p}\bigr)^{(2+p)/2p}$. 

In the case where $p=2$, the constant $\kappa$ is known explicitly:  $\kappa=\frac{\sqrt2}\pi$ (see \cite{LuPa04} and \cite{Der03}).
\end{itemize}

The article is outlined as follows. In Sections \ref{sec2} to \ref{sec5} we consider the approximation problems under the supremum norm. We start in Section \ref{sec2} by  introducting  a coding scheme which plays an important role in the sequel. 
In Section \ref{sec3}, we use the construction of Section \ref{sec2} and the self similarity of $X$ to establish a polynomial decay for $D^{(e)}(\cdot|\infty)$. In the following section, the asymptotics of the quantization error are computed. The proof relies  on a  concentration property for the entropies of ``good" coding schemes (Proposition \ref{le0811-1}).
In Section \ref{sec5}, we use the equivalence of moments in the quantization problem to establish a lower bound for  the entropy coding problem.
In the last section, we treat the case where the distortion is based on the $L^p[0,1]$-norm, i.e.\ $d(f,g)=\|f-g\|_{L^p[0,1]}$;  we introduce the distortion rate function and prove Theorem \ref{theo0303-2} with the help of Shannon's source coding Theorem.

It is convenient to use the symbols $\sim$, $\lesssim$ and $\approx$.
We write   $f \sim g$ iff
$\lim \frac fg = 1$,
while $f \lesssim g$ stands for $\limsup \frac fg \le 1$.
Finally, $f\approx g$ means
\[ 0< \liminf \frac fg \le \limsup \frac fg <\infty\ . \]

\section{The coding scheme}\label{sec2}

This section is devoted to the construction of strategies $\pi^{(n)}:\IC[0,n]\to \IL[0,n]$ which we will need later in our discussion.
The construction depends on three parameters:  $M\iN\backslash\{1\}$, $d>0$ and a strategy $\pi:\IC[0,1]\to\IL[0,1]$.

We define the maps by induction. Let $w\in \IC[0,\infty)$ and set $(w^{(n)}_t)_{t\in[0,1]}:=(w_{t+n}-w_{n})_{t\in[0,1]}$ and
$\hat w_t:=\pi(w^{(0)})(t)$ for $t\in[0,1)$. 
Assume  that $(\hat w_t)_{t\in[0,n)}$ ($n\iN$) has already been defined.
Then we choose  $\xi_n$ to be the smallest number in $\{-d+2kd/(M-1):k=0,\dots,M-1\}$  minimizing
$$
|w_{n}-(\hat w_{n-}+\xi_n)|,
$$
and extend the definition of $\hat w$ on $[n,(n+1))$ by setting
$$
\hat w_{n+t}:=\hat w_{n-}+\xi_n+\pi (w^{(n)})(t),\qquad t\in[0,1).
$$
Note that  $(\hat w_t)_{t\in[0,n)}$ depends only upon $(w_t)_{t\in[0,n)}$, so that the above construction induces strategies
$$
\pi^{(n)}:\IC[0,n]\to\IL[0,n],\ w\mapsto (\bar w^{(n)})_{t\in[0,n]},
$$
where $\bar w^{(n)}_t=\hat w_t$ for $t\in[0,n)$ and $\bar w_n= \hat w_{n-}$.
Moreover, we can write
\begin{align}\label{eq0803-4}
(\bar w_t)_{t\in[0,n]}=\pi^{(n)}(w)= \vphi_n(\pi(w^{(0)}),\dots,\pi(w^{(n-1)}), \xi_1,\dots,\xi_{n-1})
\end{align}
for an appropriate measurable function $\vphi_n:(\IL[0,n])^n\times \IR^{n-1}\to\IL[0,n]$.

The main motivation for this construction is the following property. If one has, for some $(w_t)\in\IC[0,\infty)$ and $n\iN$,
$$
\bigl\| \|w-\pi^{(n)}(w)\|_{[0,n]} \bigr\|_\infty \le \frac M{M-1}\,d 
$$
and $\|w^{(n)}-\pi(w^{(n)})\|_{[0,1]}\le d$, then
$$
|w_{n}-(\hat w_{n-}+\xi_n)|\le \frac d{M-1},
$$
whence,
\begin{align*}
\|w-\hat w\|_{[n,n+1)}&=\|w_{n}+w^{(n)}_t- (\hat
w_{n-}+\xi_n +\pi (w^{(n)})(t))\|_{[0,1)}\\
&\le  |w_{n}- (\hat
w_{n-}+\xi_n )| +
 \|w^{(n)}- \pi (w^{(n)})\|_{[0,1)}\\
&\le d/(M-1)+d=\frac{M}{M-1}\,d.
\end{align*}
In particular, if $\pi:\IC[0,1]\to\IL[0,1]$ satisfies 
$$
\bigl\|\|\fB-\pi(\fB)\|_{[0,1]}\bigr\|_\infty \le d,
$$
then for any $n\iN$,
\begin{align}\label{eq0803-3}
\bigl\|\|\fB-\pi^{(n)}(\fB)\|_{[0,n]}\bigr\|_\infty \le \frac M{M-1}\,d.
\end{align}

\section{Polynomial decay of $D^{(e)}(r|\infty)$}\label{sec3}

The objective of this section is to prove the following theorem.

\begin{theo} \label{le0803-1} There exists a constant $\kappa=\kappa(H)\in(0,\infty)$ such that
\begin{align}\label{eq1208-3}
\lim_{r\to\infty} r^H\, D^{(e)}(r|\infty)=\kappa.
\end{align}
\end{theo}

Thereafter, $\kappa=\kappa(H)$ will always denote the finite constant defined via equation (\ref{eq1208-3}).  In order to simplify notations, we abridge $\|\cdot\|=\|\cdot\|_{[0,1]}$.

\begin{remark}\label{re0207}
It was found in \cite{Der03} (see Theorem 3.5.2) that for finite moments $q\ge1$ the entropy coding error
is related to the asymptotic behavior of the  small ball function of the Gaussian measure. In particular, for fractional Brownian motion, one obtains that
$$
D^{(e)}(r|q)\approx \frac 1{r^H},\qquad r\to\infty.
$$
In order to show that $D^{(e)}(r|\infty)$ is of the order $r^{-H}$, we still need to prove an appropriate upper bound. We prove a stronger statement which will be useful later on.
\end{remark}

\begin{lemma}\label{le0802-2}
 There exist strategies $\pi^{(r)}:\IC[0,1]\to\IC[0,1]$, $r\ge0$, and probability weights $(p_w^{(r)})_{w\in \im (\pi^{(r)})}$ such that for any $q\ge1$,
\begin{align}\label{eq0314-1}
\bigl\| \|\fB-\pi^{(r)}(\fB)\| \bigr\|_\infty\le \frac1{r^H} \ \ \text{ and } \ \ \IE[(-\log p^{(r)}_{\pi^{(r)}(\fB)})^q]^{1/q} \approx r.
\end{align}
In particular, $D^{(e)}(r|\infty)\approx r^{-H}$.
\end{lemma}

The proof of the lemma is based on an asymptotic estimate for the mass concentration in randomly centered small balls, to be found in \cite{DerLif04}.
Let $\tilde X_1$ denote a fractional Brownian motion that is independent of $X$ with $\cL(X)=\cL(\tilde X_1)$. Then, for any $q\in[1,\infty)$, one has 
\begin{align}\label{eq0304-1}
\IE[(-\log \IP(\|X-\tilde X_1\|\le \eps|X))^q]^{1/q}\approx -\log \IP(\|X\|\le \eps) \approx  \eps^{-1/H}
\end{align}
as $\eps\dto0$ (see \cite{DerLif04}, Theorem 4.2 and Corollary 4.4).

\begin{proof}
For a given $\IL[0,1]$-valued sequence $(\tilde w_n)_{n\iN\cup\{\infty\}}$, we consider the following coding strategy $\pi^{(r)}(\cdot|(\tilde w_n))$:
let 
$$
T^{(r)}(w):=T^{(r)}(w|(\tilde w_n)):=\inf\{n\iN: \|w-\tilde w_n\|\le 1/r^H\},
$$
with the convention that the infimum of the empty set is $\infty$, and set 
$$
\pi^{(r)}(w):= \pi^{(r)}(w|(\tilde w_n)):=\tilde w_{T^{(r)}(w)}.
$$
Moreover, let  $(p_n)_{n\iN}$ denote the sequence of probability weights defined as
$$
p_n = \frac6{\pi^2} \,\frac1{n^2},\qquad n\iN,
$$
and set $p_\infty:=0$.

Now we let $(\tilde \fB_n)_{n\iN\cup\{\infty\}}$ denote independent FBM's that are also independent of $\fB$, and analyze the random coding strategies $\pi^{(r)}(\cdot):=\pi^{(r)}(\cdot|(\tilde \fB_n))$. With $T^{(r)}:=T^{(r)}(\fB|(\tilde \fB_n))$ we obtain
$$
\hB^{(r)}:=\pi^{(r)}(\fB)= \tilde \fB_{T^{(r)}},
$$
and
\begin{align}\label{eq1011-1}
\IE[(-\log p_{T^{(r)}})^q]^{1/q} \le 2 \IE[(\log T^{(r)})^q]^{1/q} +\log\frac{\pi^2}6.
\end{align}
Given $\fB$, the random time $T^{(r)}$ is geometrically distributed with parameter $\IP(\|\fB-\tilde \fB_1\|\le 1/ r^H|\fB)$,
and due to Lemma \ref{le041011} there  exists a universal constant $c_1=c_1(q)<\infty $ for which
$$
\IE[(\log T^{(r)})^q|\fB]^{1/q}\le c_1\,[1+ \log \IE[T^{(r)}|\fB]]=c_1 \,[1+ \log  1/\IP(\|\fB-\tilde \fB_1\|\le 1/ r^H|\fB)].
$$
Consequently,
\begin{align}
\begin{split}\label{eq1011-3}
\IE[(\log T^{(r)})^q]^{1/q}&= \IE\bigl[ \IE[( \log T^{(r)})^q|\fB]\bigr]^{1/q}\\
&\le c_1\,\IE[(1+\log 1/\IP(\|\fB-\tilde \fB_1\|\le 1/r^H|\fB))^q]^{1/q} \\
&\le c_1 \,(1+ \IE[(-\log \IP(\|\fB-\tilde \fB_1\|\le 1/r^H|\fB))^q]^{1/q}).
\end{split}
\end{align}
Due to (\ref{eq0304-1}), one has 
$$
\IE[(-\log \IP(\|\fB-\tilde \fB_1\|\le 1/r^H|\fB))^q]^{1/q} \approx  r,
$$
so that (\ref{eq1011-1}) and (\ref{eq1011-3}) imply that $\IE[(-\log p_{T^{(r)}})^q]^{1/q}\lesssim c_2 r$ for some appropriate constant $c_2<\infty$.
In particular,  for any $r\ge0$, we can find  a $\IC[0,1]$-valued  sequence $(\tilde w^{(r)})_{n\iN}$ of pairwise different elements  such that 
$$
\IE[ (-\log p_{T^{(r)}(\fB|(\tilde w^{(r)}_n))})^q]^{1/q}\le \IE[ (-\log p_{T^{(r)}})^q]^{1/q}\lesssim c_2 \,r.
$$
Now the strategies $\pi^{(r)}(\cdot|(\tilde w^{(r)}_n))$ with associated probability weights $ p^{(r)}_{\tilde w^{(r)}_n}:=p_n$  ($n\iN$)  satisfy (\ref{eq0314-1}). Moreover, $D^{(e)}(r|\infty)\approx r^{-H}$ follows since 
$$
\IH(\pi^{(r)}(X|(\tilde w_n^{(r)})))\le \IE\bigl[ -\log p_{\pi^{(r)}(X|(\tilde w_n^{(r)})}^{(r)}\bigr].
$$
\end{proof}

Let us now use the coding scheme of Section \ref{sec2} to prove

\begin{lemma}\label{le1208-2}
Let $n\iN$, $r\ge0$ and $\Delta r\ge1$. Then
\begin{align}
D^{(e)}(n(r+\Delta r)|\infty)&\le n^{-H}\,\frac{e^{\Delta r}}{e^{\Delta r}-2}\, D^{(e)}(r|\infty).
\end{align}
\end{lemma}

\begin{proof}
Fix $\eps>0$ and let $\pi :\IC[0,1]\to\IL[0,1]$ be a strategy satisfying
$$
\big\| \|\fB-\pi (\fB)\|_{[0,1]}\big\|_\infty \le (1+\eps)
D^{(e)}(r|\infty)=:d$$
and
$$
\IH(\pi (\fB))\le r.
$$
Choose $M:=\lfl e^{\Delta r} \rfl$ and let $\pi^{(n)}$ be as in Section \ref{sec2}.
Note that $\Delta r\ge1$ guarantees that $M\ge e^{\Delta r}-1\ge e^{\Delta r}
/2$, so that
$$
\bigl\| \|\fB-\pi^{(n)}(\fB)\|_{[0,n]}\bigr\|_\infty\le \frac M{M-1}\, d
\le \frac{e^{\Delta
    r}}{e^{\Delta r}-2} (1+\eps) D^{(e)}(r|\infty). 
$$
We let $(X^{(i)}_t)_{t\in[0,1]}=(X_{i+t}-X_i)_{t\in[0,1]}$ for $i=1,\dots,n$, and $(\xi_i)_{i=1,\dots,n-1}$ be as in Section \ref{sec2} for $w=X$. Observe that, due to the representation (\ref{eq0803-4}),
\begin{align}\begin{split}\label{eq0217-3}
\IH(\pi^{(n)}(X))&\le \IH(\pi (\fB^{(0)}), \dots,\pi (\fB^{(n-1)}),
\xi_1,\dots,\xi_{n-1})\\
&\le \IH(\pi (\fB^{(0)}))+\dots+\IH(\pi (\fB^{(n-1)}))+ \IH(\xi_1,\dots,\xi_{n-1})\\
&\le nr+ \log |\range (\xi_1,\dots,\xi_{n-1}) |\le  nr+n
\log M\\
&\le n(r+\Delta r).
\end{split}
\end{align}

Now let
$$\alpha_n:\ID[0,1]\to\ID[0,n], \ f\mapsto \alpha_n(f)(s)=n^H f(s/n)
$$
and consider  the strategy 
$$\tilde \pi:\IC[0,1]\to\ID[0,1], \ f\mapsto \alpha_n^{-1}\circ\pi^{(n)}\circ \alpha_n(f).
$$
Since $\alpha_n(X)$ is again a fractional Brownian motion on $[0,n]$, it follows that, a.s.
$$
\|X-\tilde\pi(X)\|_{[0,1]}= n^{-H}\, \|\alpha_n(X) -\pi^{(n)}(\alpha_n(X))\|_{[0,n]}\le (1+\eps)\, n^{-H} \frac{e^{\Delta
    r}}{e^{\Delta r}-2} D^{(e)}(r|\infty).
$$
Moreover,
$$
\IH(\tilde\pi(X))=\IH(\alpha_n^{-1}\circ \pi^{(n)}(\alpha_n(X)))=\IH(\pi^{(n)}(X))\le r.
$$
Since $\eps>0$ is arbitrary, the proof is complete.
\end{proof}

\begin{proof}[ of Theorem \ref{le0803-1}]
For $r\ge0$, $\Delta r\ge1$ and $n\iN$, Lemma \ref{le1208-2} yields
$$
D^{(e)}(n(r+\Delta r)|\infty)
\le\frac 1{n^H}\, \frac{e^{\Delta r}}{e^{\Delta r}-2}\, D^{(e)}(r|\infty).
$$
Now set $\kappa:=\liminf_{r\to\infty} r^H\, D^{(e)}(r|\infty)$ which lies in $(0,\infty)$ due to Lemma \ref{le0802-2}. Let $\eps\in(0,1/2)$ be arbitrary, and choose $r_0,\Delta r\ge1$ such that
\begin{align*}
\begin{cases}
r_0^H\, D^{(e)}(r_0|\infty) \le (1+\eps) \kappa,\\
\Delta r\le \eps r_0 \quad \text{ and}\\
e^{-\Delta r}\le \eps.
\end{cases}
\end{align*}
Then 
\begin{align*}
D^{(e)}((1+\eps)nr_0|\infty)&\le \frac 1{n^H}\,\frac 1{1-2\eps} D^{(e)}(r_0|\infty)\\
&\le \frac 1{\bigl((1+\eps) nr_0\bigr)^H}\,\frac 1{1-2\eps} (1+\eps)^{1+H} \,
\kappa
\end{align*}
and we obtain that
$$
\limsup_{n\to\infty} \bigl((1+\eps)nr_0\bigr)^H\, D^{(e)}((1+\eps)nr_0|\infty)\le \frac {(1+\eps)^{1+H}}{1-2\eps}\,  \kappa.
$$

Let now $r\ge (1+\eps)r_0$ and introduce $\bar r=\bar r(r)=\min\{(1+\eps)n r_0: n\iN, r\le (1+\eps)n
r_0\}$ as well as $\underline r=\underline r(r)=\max \{(1+\eps)n r_0: n\iN,  (1+\eps)n
r_0\le r\}$. Using the monotonicity of $D^{(e)}(r|\infty)$, we conclude that
\begin{align*}
\limsup_{r\to\infty} r^H\, D^{(e)}(r|\infty)&\le \limsup_{r\to\infty}\, \bar
  r^H\, D^{(e)}(\underline r|\infty)\\
&\le \limsup_{r\to\infty} \, (\underline r+(1+\eps)r_0)^H\, D^{(e)}(\underline r|\infty)\\
&\le \frac {(1+\eps)^{1+H}}{1-2\eps} \,\kappa.
\end{align*}
Noticing that $\eps>0$ is arbitrary finishes the proof.
\end{proof}

\section{The quantization problem}\label{sec4}

\begin{theo}\label{th1208-1}
One has for any $q\in(0,\infty)$,
$$
D^{(q)}(r|q)\sim \kappa \frac 1{r^H},\qquad r\to\infty.
$$
\end{theo}

We need some preliminary lemmas for the proof of the theorem.

\begin{lemma}\label{le0802-1}
There exist strategies $(\pi^{(r)})_{r\ge0}$ and probability weights $(p_w^{(r)})$ such that
$$
\bigl\| \|\fB-\pi^{(r)}(\fB)\| \bigr\|_\infty\le \kappa \frac1{r^H} \ \ \text{ and } \ \ -\log p^{(r)}_{\pi^{(r)}(\fB)} \lesssim r, \ \ \text{ in probability}.
$$
\end{lemma}

\begin{proof}
Let $\eps>0$ and choose $r_0\ge2$ such that 
$$
\Bigl(\frac{r_0+1}{r_0-1}\Bigr)^{1/H}\le 1+\frac\eps2
$$
By Theorem \ref{le0803-1},
$$
D^{(e)}((1+\eps/2) r|\infty)\lesssim \kappa  \frac {r_0-1}{r_0+1}\,\frac 1{r^H}
$$
In particular, there exists $r_1\ge r_0\vee \frac2\eps \log(r_0+1)$ and a map $\pi:\IC[0,1]\to \IL{[0,1]}$ such that
$$
\bigl\| \|\fB-\pi(\fB)\|_{[0,1]} \bigr\|_\infty\le \kappa  \frac {r_0-1}{r_0}\,\frac1{r_1^H}=:d \ \ \text{ and }\ \ \IH(\pi(\fB))\le (1+\eps/2)r_1.
$$
For $n\iN$, let $\pi^{(n)}$ and $\vphi_n$ be as in Section \ref{sec2} for $M=\lceil r_0\rceil$,  $d$ and $\pi$.
Then by (\ref{eq0803-3})
\begin{align}\label{eq0304-2}
\bigl\| \|\fB-\pi^{(n)}(X)\|_{[0,n]} \bigr\|_\infty\le \kappa  \frac {(r_0-1)M}{r_0(M-1)}\,\frac1{r_1^H} \le \kappa  \,\frac1{r_1^H}.
\end{align}
For $\hat w^{(0)},\dots,\hat w^{(n-1)}\in\im(\pi)$ and $k_1,\dots,k_{n-1}\in\{-d+\frac{2kd}{M-1} :k=0,\dots,M-1\}$, let $p^{(n)}$ be defined as
$$
p^{(n)}_{\vphi_n(\hat w^{(0)},\dots,\hat w^{(n-1)},k_1,\dots,k_{n-1})}=\frac {1}{M^{n-1}} \,\prod_{i=0}^{n-1} \IP(\pi(\fB)=\hat w^{(i)}).
$$
The $(p^{(n)}_w)$ define probability weights on the image of $\vphi_n$. Moreover,
$$
-\log p^{(n)}_{(\hB_t)_{t\in[0,n]}}= (n-1) \log M - \sum_{i=0}^{n-1} \log p_{\pi(\fB^{(i)})}
$$
and the ergodic theorem implies
$$
\lim_{n\to\infty} -\frac1n \,\log p^{(n)}_{(\hB_t)_{t\in[0,n]}}= \log M +\IH(\pi(\fB)),\qquad \text{a.s.}
$$
Note that $\log M +\IH(\pi(\fB))\le (1+\eps)r_1$. 

Just as in the proof of Lemma \ref{le1208-2}, we use the self similarity of $X$ to translate the strategy $\pi^{(n)}$  into a strategy for encoding $(\fB_t)_{t\in[0,1]}$.
For $n\iN$, let
$$
\alpha_n:\IL{[0,1]}\to \IL{[0,n]},\ f\mapsto (\alpha_n f)(t)={n^H}\,f(t/n)
$$
and consider $\tilde p^{(n)}_w:=p^{(n)}_{\alpha_n(w)}$ and $\tilde\pi^{(n)}(w):=\alpha^{-1}_n\circ \pi^{(n)}\circ \alpha_n (w)$.
Then
$$
-\log \tilde p^{(n)}_{\tilde \pi^{(n)}(\fB)}=-\log p^{(n)}_{\pi^{(n)}(\alpha_n(\fB))}\lesssim (1+\eps)nr_1,\qquad \text{in probability}
$$
and by (\ref{eq0304-2})
\begin{align*}
\bigl\|\|\fB-\tilde \pi^{(n)}(\fB)\|_{[0,1]} \bigr\|_\infty & =\bigl\|\|\alpha_n^{-1} ( \alpha_n(\fB)-\pi^{(n)}(\alpha_n(\fB)))\|_{[0,1]} \bigr\|_\infty\\
& = \frac 1{n^H}\, \bigl\|\| \alpha_n(\fB)-\pi^{(n)}(\alpha_n(\fB))\|_{[0,n]} \bigr\|_\infty\\
&= \frac 1{n^H}\, \bigl\|\| \fB-\pi^{(n)}(\fB)\|_{[0,n]} \bigr\|_\infty\le \kappa \frac1{(nr_1)^H}.
\end{align*}
By choosing $\bar \pi^{(r)}=\tilde \pi^{(n)}$ and $(\bar p^{(r)})=(\tilde p^{(n)})$ for $r\in((n-1)r_1,nr_1]$, one obtains a coding scheme satisfying
$$
\bigl\| \|\fB-\bar\pi^{(r)}(\fB)\| \bigr\|_\infty\le \kappa \frac1{r^H} 
$$
and
$$
 -\log \bar p^{(r)}_{\bar \pi^{(r)}(\fB)} \lesssim (1+\eps) r, \ \ \text{ in probability},
$$
so that the assertion follows by a diagonalization argument.
\end{proof}

\begin{remark} 
In the above proof, we have constructed  a high resolution coding scheme based on a strategy $\pi:\IC[0,1]\to\ID[0,1]$, using the identity $\tilde \pi_n=\alpha_n^{-1}\circ \pi^{(n)}\circ \alpha_n$.
This coding scheme leads to a coding error which is at most
\begin{align}\label{eq0218-1}
\frac{M}{M-1} \, \bigl\| \|X-\pi(X)\|_{[0,1]}\bigr\|_\infty \, n^{-H} .
\end{align}
Moreover, the ergodic theorem implies that, for large $n$, $\tilde \pi_n(X)$ lies with probability almost one in the typical  set $\{w\in\ID[0,1]: -\log \tilde p^{(n)}_w\le n (\IH(\pi(X))+ \log M+\eps)\}$, where $\eps>0$ is arbitrarily small. This set is of size $\exp\{n (\IH(\pi(X))+\log M +\eps)\}$, and will serve as a close to optimal high resolution codebook.
It remains to control the case where $\tilde \pi_n(X)$ is not in the typical set. We will do this in the proof of Theorem \ref{th1208-1} at the end of this section (see (\ref{eq0217-5})).
\end{remark}

\begin{propo} \label{le0811-1} For $q\ge1$ there exist strategies $(\pi^{(r)})_{r\ge0}$ and probability weights $(p_w^{(r)})$ such that
\begin{align}\label{eq1208-1}
\bigl\| \|\fB-\pi^{(r)}(\fB))\| \bigr\|_\infty\le \kappa \frac1{r^H} \ \ \text{ and } \ \ \lim_{r\to\infty} \frac{\IE[(-\log p^{(r)}_{\pi^{(r)}(\fB)})^q]^{1/q}}{r} =1.
\end{align}
In addition, for any $\eps>0$ one has
\begin{align}\label{eq1208-2}
\lim_{r\to\infty}\sup_{\pi, (p_w)} \IP\Bigl( -\log p_{\pi(\fB)}\le (1-\eps)r, \|\fB-\pi(\fB)\|\le \kappa\frac1{r^H}\Bigr)=0,
\end{align}
where the supremum is taken over all strategies $\pi:\IC[0,1]\to\IL[0,1]$ and over all sequences of probability weights $(p_w)$.
\end{propo}

\begin{proof}
Let $q>1$ and let $\pi^{(r)}_1$ ($r\ge0$) be a strategy and $(p^{(r,1)}_w)$ a sequence of probability weights as in Lemma \ref{le0802-1}. Moreover, let $\pi^{(r)}_2$ and $(p^{(r,2)}_w)$ ($r\ge0$) be as in  Lemma \ref{le0802-2} for $2q$.
We consider the maps $\kappa_1^{(r)}(w):=-\log p^{(r,1)}_{\pi^{(r)}_1(w)}$ and $\kappa_2^{(r)}(w):=-\log p^{(r,2)}_{\pi^{(r)}_2(w)}$, and set
$$
\pi^{(r)}(w):=\begin{cases}
\pi^{(r)}_1(w) & \text{ if } \kappa_1^{(r)}(w)\le (1+\delta)r,\\
\pi^{(r)}_2(w) & \text{ otherwise,}
\end{cases}
$$
for some fixed $\delta>0$.
Then one obtains, for  $p^{(r)}_w=\frac12 (p^{(r,1)}_w+p^{(r,2)}_w)$ and $\cT_r:=\{w\in\IC[0,1]: \kappa^{(r)}_1(w)\le (1+\delta)r\}$,
\begin{align*}
\IE [ (-\log 2p^{(r)}_{\pi^{(r)}(\fB)})^q]^{1/q} &\le \IE[1_{\cT_r}(\fB) \kappa_1^{(r)}(X)^q]^{1/q} + 
\IE[1_{\cT_r^c}(\fB) \kappa_2^{(r)}(\fB)^q]^{1/q}\\
&\le (1+\delta) r + \IP(\fB\in\cT^c_r)^{1/2q} \, \IE[\kappa_2^{(r)}(\fB)^{2q}]^{1/2q}.
\end{align*}
The definitions of $\pi^{(r)}_1$ and $\pi^{(r)}_2$ imply that $\lim_{r\to\infty}\IP(X\in\cT^c_r)=0$ and $\IE[\kappa_2^{(r)}(X)^{2q}]^{1/2q}\approx r$. Consequently, 
$$
\IE [ (-\log p^{(r)}_{\pi^{(r)}(X)})^q]^{1/q}\lesssim (1+\delta) r.
$$
Since $\delta>0$ can be chosen arbitrarily small, a diagonalization procedure leads to strategies $\tilde \pi^{(r)}$ and probability weights $(\tilde p^{(r)}_w)$ with
$$
\bigl\| \|\fB-\tilde \pi^{(r)}(\fB)\|_{[0,1]} \bigr\|_\infty \le \kappa\,\frac 1{r^H}  \ \text{ and } \ \IE[(-\log \tilde p_{\tilde \pi^{(r)}(\fB)})^q]^{1/q}\lesssim r,
$$
which proves the first assertion.

It remains to show that for arbitrary strategies $\bar\pi^{(r)}$, $r\ge0$, and probability weights $(\bar p^{(r)}_w)$:
\begin{align}\label{eq0304-3}
\lim_{r\to\infty} \IP\Bigl( -\log \bar p^{(r)}_{\bar\pi^{(r)}(\fB)}\le (1-\eps)r, \|\fB-\bar\pi^{(r)}(\fB)\|\le \kappa\frac1{r^H}\Bigr)=0.
\end{align}
Without loss of generality, we can assume  that 
\begin{align}\label{eq0304-4}
\bigl\| \|\fB-\bar\pi^{(r)}(\fB)\|_{[0,1]} \bigr\|_\infty\le \kappa\frac1{r^H}.
\end{align}
Otherwise we modify the map $\bar \pi^{(r)}$ for all $w\in\IC[0,1]$ with $\|w-\bar\pi^{(r)}(w)\|>\kappa \,r^{-H}$ in such a way that (\ref{eq0304-4}) be valid. Hereby the probability in (\ref{eq0304-3}) increases and it suffices to prove the statement for the modified strategy.
Let us consider
$$
\pi^{(r)}(w)=\begin{cases}
\bar \pi^{(r)}(w)& \text{ if } \bar p^{(r)}_{\bar \pi^{(r)}(w)} \ge \tilde p^{(r)}_{\tilde \pi^{(r)}(w)} \\
\tilde \pi^{(r)}(w) & \text{ else.} 
\end{cases}
$$
Then the probability weights $p^{(r)}:=\frac 12(\bar p^{(r)}+\tilde p^{(r)})$ satisfy
\begin{align*}
\IE[(-\log 2 p^{(r)}_{\pi(\fB)})^q]^{1/q}&\le \IE[(-\log \tilde p^{(r)}_{\tilde\pi(\fB)})^q]^{1/q} \lesssim r.
\end{align*}
Recall that
$$
\bigl\| \|\fB-\pi^{(r)}(\fB)\|_{[0,1]} \bigr\|_\infty\le \kappa\frac1{r^H},
$$
hence by Theorem \ref{le0803-1}, one has $\IE[ -\log p^{(r)}_{\pi^{(r)}(\fB)}]\ge \IH(\pi^{(r)}(X))\gtrsim r$. Lemma \ref{le0809-1} thus implies that
$$
-\log p^{(r)}_{\pi^{(r)}(\fB)}\sim r,\quad \text{ in probability.}
$$
In particular, 
$$
-\log \bar p^{(r)}_{\bar \pi^{(r)}(\fB)}\ge -\log 2 p^{(r)}_{\pi^{(r)}(\fB)}\gtrsim r,\quad\text{ in probability,}
$$
which implies (\ref{eq0304-3}).
\end{proof}

\begin{proof}[ of Theorem \ref{th1208-1}]  We start by proving the lower bound.
Fix $q>0$,
let $\cC_r$, $r\ge0$, denote arbitrary codebooks of size $e^{r}$, and let $\pi^{(r)}:\IC[0,1]\to\cC_r$ denote  arbitrary strategies. Moreover, let $(p^{(r)}_w)$ be the sequence of probability weights defined as $p^{(r)}_w=1/|\cC_r|$, $w\in\cC_r$.
Then $-\log p^{(r)}_{\pi^{(r)}(X)}\le r$ a.s., and the above lemma implies that for any $\eps\in(0,1)$,
$$
\lim_{r\to\infty} \IP\Bigl(\|X-\pi^{(r)}(X)\|\le \kappa \frac{(1-\eps)^H}{r^H}\Bigr)=0.
$$
Therefore,
$$
\IE[\|X-\pi^{(r)}(X)\|^q]^{1/q}\ge \kappa \frac{(1-\eps)^H}{r^H}\,\IP\Bigl(\|X-\pi^{(r)}(X)\|\ge \kappa \frac{(1-\eps)^H}{r^H}\Bigr)^{1/q} \sim \kappa \frac{(1-\eps)^H}{r^H},
$$
which proves the lower bound.

It remains to show that $D^{(q)}(r,q)\lesssim \kappa/r^H$. 
By Lemma \ref{le0802-1}, there exist strategies $\pi^{(r)}$ and probability weights $(p^{(r)}_w)$ such that 
$$
\bigl\| \|\fB-\pi^{(r)}(\fB)\| \bigr\|_\infty \le \kappa \frac 1{r^H}\ \text{ and } \ -\log p_{\pi^{(r)}(\fB)} \lesssim r, \quad \text{in probability}.
$$
Furthermore, due to Theorem 4.1 in \cite{DFMS03}, there exist  codebooks $\bar \cC_r$ of size $e^r$ with
$$
\IE[\min_{\hat w\in\bar \cC_r} \|\fB-\hat w\|^{2q}]^{1/{2q}} \approx \frac 1{r^H}.
$$
We consider the codebook $\cC_r:=\bar \cC_r\cup\{\hat w: -\log p^{(r)}_{\hat w}\le (1+\eps/2)r\}$. Clearly, $\cC_r$ contains at most $e^r+e^{(1+\eps/2)r}$ elements. Moreover,
\begin{align}\begin{split}\label{eq0217-5}
\IE[\min_{\hat w\in\cC_r} \|\fB-\hat w\|^q]^{1/{q}}&\le \IE[ 1_{\cC_r}(\pi^{(r)}(\fB))\,  (\kappa \frac 1{r^H})^q]^{1/q}\\
& \ \ \ +
\IE[ 1_{\cC_r^c}(\pi^{(r)}(\fB)) \min_{\hat w\in\bar \cC_r} \|\fB-\hat w\|^q]^{1/q}\\
&\le \kappa \frac 1{r^H}+ \IP(\pi^{(r)}(\fB)\not\in \cC_r)^{1/2q} \, \IE[\min_{\hat w\in\bar \cC_r} \|\fB-\hat w\|^{2q}]^{1/2q}.
\end{split}
\end{align}
Since $\lim_{r\to\infty} \IP(\pi^{(r)}(\fB)\not\in \cC_r)=0$ and the succeeding expectation is of order $\cO(1/r^H)$, the second summand is of order $o(1/r^H)$. Therefore, for $r\ge 2/\eps$
$$
D^{(q)}((1+\eps)r|q) \le \IE[\min_{\hat w\in\cC_r} \|\fB-\hat w\|^q]^{1/{q}}\lesssim \kappa \frac 1{r^H}.
$$
By switching from $r$ to $\tilde r=(1+\eps)r$, we obtain
$$
D^{(q)}(\tilde r|q) \lesssim \kappa\, (1+\eps)^H \,\frac 1{\tilde r^H}.
$$
Since $\eps>0$ was arbitrary, the proof is complete.
\end{proof}

\section{Implications of the equivalence of  moments}\label{sec5}

In this section we complement Theorem \ref{th1208-1} by

\begin{theo}\label{theo0209-1} For arbitrary $q\in(0,\infty]$, one has
$$
D^{(e)}(r|q)\sim \kappa\frac1{r^H}.
$$
\end{theo}

The proof of this theorem is based on the following general principle: \emph{if the asymptotic quantization error coincides for two different moments $q_1<q_2$, then all moments $q\le q_2$ lead to the same asymptotic quantization error and the entropy coding problem coincides with the quantization problem for all moments $q\le q_2$.}

Let us prove this relationship in a general setting.
$E$ and  $\hat E$ denoting arbitrary measurable spaces and
$d:E\times \hat E\to[0,\infty)$ a measurable function, the quantization error for a general $E$-valued r.v.\ $X$ under the distortion $d$ is defined as
$$
D^{(q)}(r|q)=\inf_{\cC\subset E} \IE[\min_{\hat x\in\cC} d(X,\hat x)^q]^{1/q},
$$
where the infimum is taken over all codebooks $\cC\subset \hat E$ with $|\cC|\le e^r$. In order to simplify notations, we abridge 
$$
d(x,A)=\inf_{y\in A} d(x,y),\qquad x\in E,\ A\subset \hat E.
$$ 
Analogously, we denote the entropy coding error by
$$
D^{(e)}(r|q)=\inf_{\hat X} \IE[d(X,\hat X)^q]^{1/q} ,
$$
where the infimum is taken over all discrete $\hat E$-valued r.v.\ $\hat X$ with $\IH(\hat X)\le r$.

Then Theorem \ref{theo0209-1} is a consequence of Theorem \ref{th1208-1} and the following theorem.

\begin{theo}\label{theo1201-1} Assume that $f:[0,\infty)\to\IR_+$ is a decreasing, convex function satisfying
\begin{align}\label{eq1201-1}
\limsup_{r\to\infty}\frac{-r \frac{\partial^+}{\partial r} f(r)}{f(r)}<\infty,
\end{align}
and suppose that, for some $0<q_1<q_2$,
$$
D^{(q)}(r+\log 2|q_1)\sim D^{(q)}(r|q_2) \gtrsim f(r).
$$ 
Then for any $q>0$,
$$
D^{(e)}(r|q)\gtrsim f(r).
$$
\end{theo}

We need some technical lemmas.

\begin{lemma}\label{le1201-2}
Let $0<q_1<q_2$ and $f:[0,\infty) \to\IR_+$. If 
$$
D^{(q)}(r+\log2|q_1)\sim D^{(q)}(r|q_2)\sim f(r),
$$
then for any $\eps>0$,
$$
\lim_{r\to\infty} \suptwo{\cC\subset E:}{ |\cC|\le e^r} \IP(d(X,\cC) \le (1-\eps) f(r))=0.
$$
\end{lemma}

\begin{proof}
For $r\ge0$, let $\cC_r^*$ denote codebooks of size $e^r$ with
\begin{align}\label{eq0908-2}
\IE[d(X,\cC^*_r)^{q_2}]^{1/q_2}\sim f(r).
\end{align}
Now let $\cC_r$ denote arbitrary codebooks of size $e^r$, and consider the codebooks $\bar\cC_r:=\cC^*_r\cup  \cC_r$.
Using (\ref{eq0908-2}) and the inequality  $q_1\le q_2$, it follows that
$$
f(r)\gtrsim \IE[d(X,\bar \cC_r)^{q_2}]^{1/q_2}\ge \IE[d(X,\bar \cC_r)^{q_1}]^{1/q_1} \ge D^{(q)}(r+\log 2|q_1)\sim f(r).
$$
Hence, Lemma \ref{le0809-1} implies that
$$
d(X,\bar \cC_r) \sim f(r),\quad \text{in probability},
$$
so that in particular,
$$
d(X,\cC_r) \gtrsim f(r),\quad \text{in probability}.
$$
\end{proof}

\begin{lemma}\label{le1201-1} Assume that $f:[0,\infty)\to\IR_+$ is a decreasing, convex function satisfying (\ref{eq1201-1}) and
$$
\lim_{r\to\infty} \suptwo{\cC\subset E:}{ |\cC|\le e^r} \IP(d(X,\cC) \le f(r))=0.
$$
Then for any $q>0$,
$$
D^{(e)}(r|q)\gtrsim f(r).  
$$
\end{lemma}

\begin{proof}
The result is a consequence of the technical Lemma \ref{le0811-2}.
Consider the family $\cF$ consisting of all random vectors
$$
(A,B)=(d(X,\hat X)^q, -\log p_{\hat X}), 
$$
where $\hat X$ is an arbitrary discrete $E$-valued r.v.\ and $(p_w)$ is an arbitrary sequence of probability weights on the range of $\hat X$. 
Let $\tilde f(r)=  f(r)^q$, $r\ge0$.
Then for any choice of $\hat X$ and $(p_w)$ and an arbitrary $r\ge0$, the set $\cC:=\{w\in E: -\log p_w\le r\}$ contains at most $e^r$ elements. Consequently,
$$
\IP( d(X,\hat X)^q\le \tilde f(r), -\log p_{\hat X}\le r)= \IP( d(X,\hat X)\le f(r), \hat X\in\cC)\le \IP( d(X,\cC)\le f(r)).
$$
By assumption the right hand side converges  to $0$ as $r\to\infty$ ,independently of the choice of $\hat X$ and $(p_w)$. 
Since $\tilde f$ satisfies condition (\ref{condition}), Lemma \ref{le0811-2}  implies that
\begin{align*}
D^{(e)}(r|q)= \inf_{\hat X: \IH(\hat X)\le r} \IE[d(X,\hat X)^q]^{1/q} =\inf_{A\in\cF_r} \IE [A]^{1/q}\gtrsim \tilde f(r)^{1/q}=f(r), 
\end{align*}
where $\cF_r=\{A:(A,B)\in\cF,\ \IE B\le r\}$.
\end{proof}

Theorem \ref{theo1201-1} is now an immediate consequence of Lemma \ref{le1201-2} and Lemma \ref{le1201-1}.

\section{Coding with repect to the $L^p[0,1]$-norm distortion}\label{sec6}

In this section, $p\in[1,\infty)$ is fixed. In contrast to the previous sections, we consider entropy coding and quantization of $X$ in $L^p[0,1]$, i.e.\ $\hat E=L^p[0,1]$ and $d(f,g)=\|f-g\|_{L^p[0,1]}$.
In order to treat these approximation problems, we need to introduce Shannon's distortion rate function. It is defined as
$$
D(r|q)=\inf \bigl\| \|X-\hat X\|_{L^p[0,1]} \bigr\|_q,
$$
where the infimum is taken over all $\hat E$-valued r.v.'s $\hat X$ satisfying the mutual information constraint $I(X;\hat X)\le r$. Here and elsewhere $I$ denotes the \emph{Shannon mutual information}, defined as
$$
I(X;\hat X) =\begin{cases}
\int \log \frac{d\IP_{X,\hat X}}{d\IP_X\otimes \IP_{\hat X}} \, d\IP_{X,\hat X}  & \text{if } \IP_{X,\hat X}\ll \IP_X\otimes \IP_{\hat X}\\
\infty &\text{else}.
\end{cases}
$$
The objective of this section is to prove

\begin{theo}\label{theo0207-1}
The following limit exists
\begin{align}\label{eq0207-1}
\kappa_p=\kappa_p(H)=\lim_{r\to\infty} r^H\,D(r|p)\in(0,\infty),
\end{align}
and for any $q>0$, one has
\begin{align}\label{eq0209-1}
D^{(q)}(r|q)\sim D^{(e)}(r|q)\sim \kappa_p\,\frac1{r^H}.
\end{align}
\end{theo}

We will first prove that statement (\ref{eq0209-1}) is valid for
$$
\kappa_p:=\liminf_{r\to\infty} r^H\,D(r|p).
$$
Since  $D(r|p)$ is dominated by $D^{(q)}(r|p)$, the existence of the limit in (\ref{eq0207-1}) then follows immediately. 
Due to Theorem 1.2 in \cite{Der04}, the distortion rate function $D(\cdot|p)$ has the same weak asymptotics as $D^{(q)}(\cdot|p)$. In particular, $D(r|p)\approx r^{-H}$ and $\kappa_p$ lies in $(0,\infty)$.

We proceed as follows: decomposing $X$ into the two processes
$$
X^{(1)}=(X_t-X_{\lfl t\rfl})_{t\ge0} \ \ \text{ and } \ \ X^{(2)}=(X_{\lfl t\rfl})_{t\ge0},
$$
we consider the coding problem for $X^{(1)}$ and $X^{(2)}$ in $L^p[0,n]$ ($n\iN$ being large).
We control the coding complexity of the first term via Shannon's Source Coding Theorem (SCT) and use a limit argument in order to show that the coding complexity of $X^{(2)}$ is asymptotically negligible. 
We recall the SCT in a form which is appropriate for our discussion;
for $n\iN$, let
$$
d_{p}(f,g)=\Bigl(\int_0^1 |f(t)-g(t)|^p \,{dt} \Bigr)^{1/p}
$$
and
$$
d_{n,p}(f,g)=\Bigl(\int_0^n |f(t)-g(t)|^p \,\frac{dt}n \Bigr)^{1/p}.
$$
Then $\tilde d_n(f,g)=d_{n,p}(f,g)^p$, $n\iN$, is a single letter distortion measure, when  interpreting the function $f|_{[0,n)}$ as  the concatenation of the ``letters" $f^{(0)},\dots, f^{(n-1)}$, where $f^{(i)}=(f(i+t))_{t\in[0,1)}$. Analogously, the process $X^{(1)}$ corresponds to the letters  $X^{(1,i)}:=(X_{i+t})_{t\in[0,1)}$, $i\iN_0$.  Since $(X^{(1,i)})_{i\iN_0}$ is an ergodic stationary $\IC[0,1)$-valued process, the SCT implies that
for fixed $r\ge0$ and $\eps>0$ there exist codebooks $\cC_n\subset L^p[0,n]$, $n\iN$, with at most $\exp\{(1+\eps)nr\}$ elements such that
\begin{align}\label{eq0302-1}
\lim_{n\to\infty} \IP(\tilde d_n(X^{(1)},\cC_n) \le (1+\eps) D(r|p)^p)=1.
\end{align}
A proof of this statement can be carried out by using the asymptotic equipartition property as  stated in \cite{DeKo02} (Theorem 1). The proof is standard and therefore omitted. For further details concerning the distortion rate function one can consult \cite{CoTho91} or \cite{DeKo02}.  

First we prove a lemma which will later be used to control the coding complexity of \nolinebreak[4] $X^{(2)}$.

\begin{lemma}\label{le0302-1} Let $(Z_i)_{i\iN}$ be an ergodic stationary sequence of real-valued r.v.'s and let $S_n=\sum_{i=1}^n Z_i$, $n\iN_0$.
Then there exist codebooks $\cC_n\subset \IR^n$ of size $\exp\{n\IE[ \log (|Z_1|/2\eps+2)]+ nc\}$ satisfying
$$
\lim_{n\to\infty} \IP\bigl( \min_{\hat s\in\cC} \|S_1^n-\hat s\|_{l^n_\infty} )\le \eps\bigr)=1,
$$
where $S_1^n$ denotes $(S_i)_{i=1,\dots,n}$,  $c$ is a universal constant and $\|\cdot\|_{l^n_\infty}$ denotes the maximum norm on $\IR^n$.
\end{lemma}

\begin{proof}
Let $c>0$ be such that 
$(p_n)_{n\iZ}$ defined through
$$
p_n = e^{-c}\,  \frac 1{(|n|+1)^2}
$$
is a sequence of probability weights. For a given sequence $(s_n)_{n\iN}$, we define a reconstruction $(\hat s_n)$ recursively. The construction depends on a parameter $\eps>0$. Let $\hat s_0=0$ and suppose that $\hat s_0^n=(\hat s_i)_{i=0,\dots,n}$ is already defined. Then we choose a $\xi_{n+1}\in 2\eps \IR$ minimizing the distance
$$
|s_{n+1} - (\hat s_{n}+\xi_{n+1})|
$$
and set $\hat s_{n+1}:=\hat s_{n}+\xi_{n+1}$.
This defines maps $\pi_n:\IR^n \to \IR^n, s_1^n \mapsto \pi_n(s_1^n):=\hat s_1^n$.
We equip the range of $\pi_n$ with a sequence of probability weights via
$$
p^{(n)}_{\hat s_1^n} =\prod_{i=1}^n p_{\xi_i/2\eps}.
$$
Then
\begin{align*}
-\log p^{(n)}_{\hat s_1^n}\le 2 \sum_{i=1}^n \log (|\xi_i|/2\eps+1)+ n\,c.
\end{align*}
Now consider $\pi_n(S_1^n)$. Let $\xi_n=\xi_n((S_i))$ be as above when replacing the deterministic argument $(s_n)$ by $(S_n)$.
Then 
$$
|\xi_n-Z_n|= | \hat S_n-\hat S_{n-1}- S_n+S_{n-1}|\le 2\eps 
$$
and, hence,
$|\xi_n|\le |Z_n|+2\eps$. Consequently,
\begin{align*}
-\frac1n \log p^{(n)}_{\hat S_1^n}\le 2 \frac1n \sum_{i=1}^n \log (|Z_i|/2\eps+2)+ c\to 2\IE[ \log (|Z_1|/2\eps+2)]+ c,
\end{align*}
where the convergence follows due to the ergodicity of $(Z_n)$.
Therefore the codebooks
$$
\cC_n:=\bigl\{\hat s_1^n\iR^n: -\frac1n \log p^{(n)}_{\hat s_1^n}\le 2 \IE[ \log (|Z_1|/2\eps+2)]+ 2c\bigr\} 
$$
satisfy the required assertion.
\end{proof}

We now use the SCT combined with the previous lemma to construct codebooks that guarantee almost optimal reconstructions with a high probability.

\begin{lemma}
For any $\eps>0$ there exist codebooks $\cC_r$, $r\ge0$, of size $e^r$ such that
$$
\lim_{r\to\infty} \IP(d_p(X,\cC_r)\le (1+\eps) \kappa_p r^{-H})=1.
$$
\end{lemma}

\begin{proof} Let $\eps>0$ be arbitrary and $c$ be as in Lemma \ref{le0302-1}.
We fix $r_0\ge \bigl( \frac{4\eps\kappa_p}{\IE|X_1|}\bigr)^{1/H}$ such that
\begin{align}\label{eq1221-1}
\eps\kappa_p r^{-H}\ge e^{-\eps r+c+\log\IE |X_1|}
\end{align}
for all $r\ge r_0$. Then choose $r_1\ge r_0$ with
$$
D(r_1|p)\le (1+\eps) \kappa_p r_1^{-H}.
$$
We decompose $X$ into the two processes
$$
X^{(1)}_t=X_t-X_{\lfl t\rfl} \ \ \text{ and } \ \ X^{(2)}_t=X_{\lfl t\rfl}.
$$
Due to the SCT (\ref{eq0302-1}), there exist codebooks $\cC^{(1)}_n\subset L^p[0,n]$ of size $\exp\{ (1+\eps)nr_1\}$ satisfying
$$
\lim_{n\to\infty} \IP( d_{n,p}(X^{(1)},\cC^{(1)}_n)^p\le (1+2\eps)^p  \kappa_p^p r_1^{-pH})=1.
$$ 
We apply Lemma \ref{le0302-1} for $\eps':=\eps\kappa_p r_1^{-H}$. Note that
$$
\IE \log \Bigr(\frac {|X_1|}{2\eps'}+2\Bigr)+c \le \log \Bigr(\frac {\IE |X_1|}{2\eps'}+2\Bigr)+c
$$
Since $r_1^H \ge  \frac{4\eps\kappa_p}{\IE|X_1|}$, it follows that $\frac {\IE|X_1|}{2\eps'}= \frac {r_1^H\,\IE|X_1|}{2\eps\kappa_p}\ge2$, so that
\begin{align*}
\IE \log \Bigr(\frac {|X_1|}{2\eps'}+2\Bigr)+c&\le \log \Bigr(\frac {\IE |X_1|}{\eps'}\Bigr)+c\\
&= -\log (\eps\kappa_p r_1^{-H}) +c+ \log \IE |X_1|\le\eps r,
\end{align*}
due to (\ref{eq1221-1}).
Hence, there exist codebooks $\cC^{(2)}_n\subset L^p[0,n]$ of size $\exp\{\eps n r_1\}$ with
$$
\lim_{n\to\infty} \IP\Bigl(d_{n,p}(X^{(2)},\cC^{(2)}_n)\le \eps \kappa_p \frac 1{r_1^H}\Bigr)=1.
$$ 
Let now $\tilde\cC_n:=\cC_n^{(1)}+\cC_n^{(2)}$ denote the Minkowski sum of the sets $\cC_n^{(1)}$ and $\cC_n^{(2)}$. Then $|\tilde \cC_n|\le \exp\{(1+2\eps)nr_1\}$, and one has
\begin{align*}
\IP(d_{n,p}(X,\tilde \cC_n) \le (1+3\eps)  \kappa_p r_1^{-H}) \ge \IP(  & d_{n,p}(X^{(1)},\cC_n^{(1)}) \le (1+2\eps)  \kappa_p r_1^{-H} \text{ and}\\
&  d_{n,p}(X^{(2)},\cC_n^{(2)}) \le \eps  \kappa_p r_1^{-H})\to 1.
\end{align*}
Consider the isometric isomorphism
$$
\beta_n:L^p[0,1]\to (L^p[0,n],d_{n,p}),\  f\mapsto f(nt),
$$
and the codebooks $\cC_n\subset L^p[0,1]$ given by
$$
\cC_n=\{n^{-H} \beta_n^{-1} (\hat w): \hat w\in\tilde\cC_n\}
$$
Then $\tilde X^{(n)}= n^{-H} \beta^{-1}_n (X)$ is a fractional Brownian motion and one has
\begin{align*}
d_p(\tilde X^{(n)},\cC_n)= d_{n,p}(\beta_n(\tilde X^{(n)}), \beta_n(\cC_n))= n^{-H}  d_{n,p}(X,\tilde \cC_n).
\end{align*}
Hence, the codebooks $\cC_n$ are of size $\exp\{(1+2\eps)nr_1\}$ and satisfy
$$
\IP(d_p(X,\cC_n)\le (1+3\eps)  \kappa_p (n r_1)^{-H}))= \IP(   d_{n,p}(X,\tilde \cC_n)\le (1+3\eps)  \kappa_p  r_1^{-H})\to 0
$$
as $n\to\infty$. Now the general statement follows by an interpolation argument similar to that used at the end of the proof of Theorem \ref{le0803-1}.
\end{proof}

\begin{proof}[ of Theorem \ref{theo0207-1}]
Let $q\ge1$ be arbitrary, let $\cC^{(1)}_r$ be  as in the above lemma for some fixed $\eps>0$. Moreover, we let $\cC^{(2)}_r$ denote codebooks of size $e^r$ with
$$
\IE[ d_p(X,\cC^{(2)}_r)^{2q}]^{1/(2q)}\approx \frac1{r^H}.
$$
Then the codebooks $\cC_r:=\cC_r^{(1)}\cup \cC_r^{(2)}$ contain at most $2e^r$ elements and satisfy, in analogy to the proof of Theorem \ref{th1208-1} (see (\ref{eq0217-5})),
$$
\IE[d_p(X,\cC_r)^q]^{1/q}\lesssim (1+\eps) \kappa_p \frac1{r^H},\qquad r\to\infty.
$$
Since $\eps>0$ is arbitrary, it follows that
$$
D^{(q)}(r|q)\lesssim \kappa_p \frac1{r^H}.
$$
For $q\ge p$ the quantization error is greater than the distortion rate function $D(r|p)$, so that the former inequality extends to
$$
\lim_{r\to\infty} r^H\,D^{(q)}(r|q)=\kappa_p.
$$
In particular, we obtain the asymptotic equivalence of all moments $q_1,q_2$ greater or equal to $p$.
Next, an application of Theorem  \ref{theo1201-1} with $d(f,g)= d_p(f,g)^q$ implies that for any $q>0$,
$$
D^{(e)}(r|q)\gtrsim \kappa_p \frac 1{r^H},
$$
which establishes the assertion.
\end{proof}

\begin{appendix}

\section*{Appendix}

\setcounter{section}{1}
\setcounter{theo}{0}

\begin{lemma} \label{le0809-1} For $r\ge0$, let $A_r$  denote $[0,\infty)$-valued r.v.'s. If one has, for $0<q_1<q_2$ and some function $f:[0,\infty)\to\IR_+$,
\begin{align}\label{eq0908-1}
\IE[A_r^{q_1}]^{1/q_1} \sim \IE[A_r^{q_2}]^{1/q_2}\sim f(r), 
\end{align}
then 
$$
A_r \sim f(r), \text{ in probability}.
$$
\end{lemma}

\begin{proof}
Consider 
$$
\tilde A_r:=A_r^{q_1}/\IE[A_r^{q_1}],
$$
and $\tilde q_2=q_2/q_1$. Then (\ref{eq0908-1}) implies that
$$
\IE[\tilde A_r^{\tilde q_2}]^{1/\tilde q_2}\sim \IE[\tilde A_r]=1 
$$
Denoting  $\Delta \tilde A_r:=\tilde A_r-1$ and $g(x):=x^{\tilde q_2}$, we obtain
\begin{align*}
\IE[\tilde A_r ^{\tilde q_2}]&=\IE[1+\Delta \tilde A_r g'(1) + g(1+\Delta\tilde A_r)-(1+\Delta \tilde A_r \, g'(1))]\\
&= 1 + \IE[ g(\tilde A_r)- (1+\Delta \tilde A_r \,g'(1))]
\end{align*}
Due to the strict convexity of $g$, for  arbitrary $\eps>0$ there exists $\delta>0$ such that
$$
g(x+1)\ge 1+x g'(1)+\delta,  \ \text{ for } \  x\in[-1,1-\eps]\cup [1+\eps,\infty).
$$
Consequently,
$$
\IE[\tilde A_r ^{\tilde q_2}]\ge 1+ \delta \,\IP(|\Delta\tilde A_r|\ge \eps).
$$
Since $\lim_{r\to\infty}\IE[\tilde A_r ^{\tilde q_2}]=1$, it follows that $\lim_{r\to\infty}\IP(|\Delta\tilde A_r|\ge \eps)=0$.
Hence, 
$$
A_r = \IE[A_r^{q_1}]^{1/q_1}\, \tilde A_r^{1/q_1}\sim \IE[A_r^{q_1}]^{1/q_1} \sim f(r),\quad \text{in probability}.
$$
\end{proof}

\begin{lemma} \label{le041011} Let $q\ge1$. There exists a constant $c=c(q)<\infty$ such that
for all $[1,\infty)$-valued r.v.'s $Z$ one has
$$
\IE[ (\log Z)^q]^{1/q}\le c\,[1+ \log \IE [Z]].
$$
\end{lemma}

\begin{proof} Using elementary analysis, there exists a positive constant $c_1=c_1(q)<\infty$  such that $\psi(x):= (\log x)^q+c_1 \,\log x$, $x\in[1,\infty)$, is concave. For any $[1,\infty)$-valued r.v.\ $Z$, Jensen's inequality then yields
\begin{align*}
\IE[ (\log Z)^q]^{1/q} & \le \IE[\psi(Z)]^{1/q}\le \psi(\IE[Z])^{1/q}\\
& \le \log\IE[Z] + c_1^{1/q} (\log \IE[Z])^{1/q}\le c\, [1+ \log\IE[Z]],
\end{align*}
where $c=c(q)<\infty$ is an appropriate universal constant.
\end{proof}

\begin{lemma}\label{le0811-2}
Let $f:[0,\infty) \to\IR_+$ be a decreasing, convex function satisfying $\lim_{r\to\infty} f(r)=0$ and
\begin{equation}\label{condition}
\limsup_{r\to\infty} \frac {-r \,\frac{\partial^+}{\partial r} f(r)}{f(r)}<\infty,
\end{equation}
and $\cF$ be a family of $[0,\infty]^2$-valued random variables for which
\begin{align}\label{eq0811-4}
\lim_{r\to\infty} \sup_{(A,B)\in\cF} \IP(A\le f(r), B\le r)=0.
\end{align}
Then the sets of random variables $\cF_r$ defined for $r\ge0$ through 
$$
\cF_r:=\{A: (A,B)\in\cF,\ \IE B\le r\}
$$
satisfy 
$$
\inf_{A\in\cF_r}  \IE A \gtrsim f(r)
$$
as $r\to\infty$.
\end{lemma}

\begin{proof}
Fix $R>0$, positive integers $I$ and $N$, and define $\lambda := - \frac{\partial^+}{\partial r} f(R)$,
$$
r_i:= \frac{i+N}N  R, \qquad i=-N,-N+1,\dots. 
$$
For $(A,B) \in \cF_R$, we define
$$
\cT_{A,B}:=\{\nexists i\in\{-N+1,\dots, I\} \text{ such that } A \le f(r_i) \text { and } B \le r_i\}.
$$
Then we have
\begin{eqnarray*}
\IE\bigl[A+\lambda B\bigr] &\ge& \sum_{i=-N}^{I-1} \IE \bigl[1_{\cT_{A,B}} 1_{[r_i,r_{i+1})}(B)(A+\lam r_i)\bigr]\\
&\ge&  \sum_{i=-N}^{I-1} \IE \bigl[1_{\cT_{A,B}} 1_{[r_i,r_{i+1})}(B)(f(r_{i+1})+\lam r_i)\bigr]\\
&=&  \sum_{i=-N}^{I-1} \IE \Bigl[1_{\cT_{A,B}} 1_{[r_i,r_{i+1})}(B)(f(r_{i+1})+\lam r_{i+1}-\lam \frac{R}{N})\Bigr]\\ 
&\ge&  \sum_{i=-N}^{I-1} \IE \Bigl[1_{\cT_{A,B}} 1_{[r_i,r_{i+1})}(B)(f(R)+\lam R -\lam \frac{R}{N})\Bigr],\\ 
\end{eqnarray*}
where the last inequality follows from the fact that
$$
f(R)+\lam R = \inf_{r \ge 0} \left[ f(r) + \lam r \right]
$$ 
by the definition of $\lam$ and the convexity of $f$.
Now,  fix $\eps >0$ and pick $N \ge 1/\eps$, $I \ge 
2N/\eps$ and $R_0$ so large that 
$$
\IP  (\cT_{A,B}) \ge 1 - \frac{\eps}{2} \mbox{ for all } R \ge R_0 \mbox{ and all } (A,B) \in \cF_R.
$$
Using Chebychev's inequality, we then obtain for $R \ge R_0$,
\begin{eqnarray*}
\IE [ A+\lam B] &\ge& (1-\eps )(f(R)+\lam R) \left( 1-\IP \left( \cT ^c \right) - \IP \left( B \ge R \frac{I}{N} \right) \right)\\
&\ge&  (1-\eps )(f(R)+\lam R) \left( 1 - \frac{\eps}{2} - \frac{\eps}{2} \right) . 
\end{eqnarray*}
Hence,
$$
\lam R + \IE A \ge (1-\eps )^2 \left( f(R) + \lam R \right)
$$
and therefore
$$
\IE A \ge (1-\eps )^2 f(R)   + \lam R \left(  (1-\eps )^2-1 \right).
$$
Using the definition of $\lam$ and  (\ref{condition}), as well as the fact that $\eps>0$ is arbitrary,
the conclusion follows.
\end{proof}

\end{appendix}

\bibliography{biblio}
\bibliographystyle{abbrv}

\end{document}